\documentclass[letterpaper, 10 pt, conference]{ieeeconf}  %

\IEEEoverridecommandlockouts                              %
\usepackage{multicol}
\usepackage[pageanchor=true,plainpages=false, pdfpagelabels, bookmarks,bookmarksnumbered,hidelinks]{hyperref}
\usepackage{pdfsync}		%
\usepackage{amsmath}		%
\usepackage{amssymb}		%
\usepackage{bm}             %
\usepackage{graphicx}		%
\usepackage[footnotesize]{caption}
\usepackage[listofformat=subsimple,
labelformat=simple]{subfig}
\usepackage{psfrag}   	%
\usepackage{color}			%
\usepackage{balance}    %
\usepackage[sort]{cite}
\usepackage{mathnotation} %
\usepackage{nopageno}

\newcommand{\subparagraph}{}
\usepackage{titlesec}
\titlespacing\section{0pt}{10pt plus 0pt minus 2pt}{0pt plus 2pt minus 2pt}
\titlespacing\subsection{0pt}{10pt plus 0pt minus 2pt}{0pt plus 2pt minus 2pt}
\titlespacing\subsubsection{0pt}{10pt plus 0pt minus 2pt}{0pt plus 2pt minus 2pt}

\usepackage{marginnote}

\title{\LARGE \bf
Soap-bubble Optimization of Gaits
}

\author{Suresh Ramasamy and Ross L. Hatton \\ Oregon State University \\ \{ramasams,hattonr\}@oregonstate.edu%
}

\begin{document}

\maketitle
\thispagestyle{empty}
\pagestyle{plain}

\begin{abstract}
In this paper, we present a geometric variational algorithm for optimizing the gaits of kinematic locomoting systems. The dynamics of this algorithm are analogous to the physics of a soap bubble, with the system's Lie bracket supplying an ``inflation pressure'' that is balanced by a ``surface tension'' term derived from a Riemannian metric on the system's shape space. We demonstrate this optimizer on a variety of system geometries (including Purcell's swimmer) and for optimization criteria that include maximizing displacement and efficiency of motion for both translation and turning motions.

\end{abstract}

\section{Introduction}

Gait optimization algorithms for locomoting systems must contend with a number of nonlinear effects, including:
\bi
\ie Shape-dependent system dynamics (and thus a control structure that depends on the state of the controller and not just the plant); and
\ie History-dependent input-output mappings (e.g. if one part of the gait rotates the system, future translational sections will propel the system in a new direction).
\ei
The geometric mechanics community has long recognized that Lie brackets of the system dynamics can be used to provide an understanding of these nonlinearities, and thereby to identify shape oscillations that produce useful net displacements through the world~\cite{Murray:1993,Morgansen:2007}. For kinematic locomoting systems with two shape variables (a class which includes a number of informative systems in viscous, granular, inertial, and fluid environments~\cite{walsh95, ostrowski98a,Shammas:2007,Avron:2008,Hatton:2013TRO:Swimming,Melli:2006,Kanso:2009sw,Hatton:2013PRL}), the Lie bracket can be visualized as a set of ``height functions'' over the shape space, whose topographies identify good placements for gait trajectories. For example, gaits that move a system forward encircle positive peaks (or negative valleys) in the $x$ component of the Lie bracket, whereas gaits that turn the system with minimal translation enclose sign-definite regions in the $\theta$ component and balanced positive and negative regions in the $x$ and $y$ components.

\begin{figure}%
\centering
\includegraphics[width=.45\textwidth]{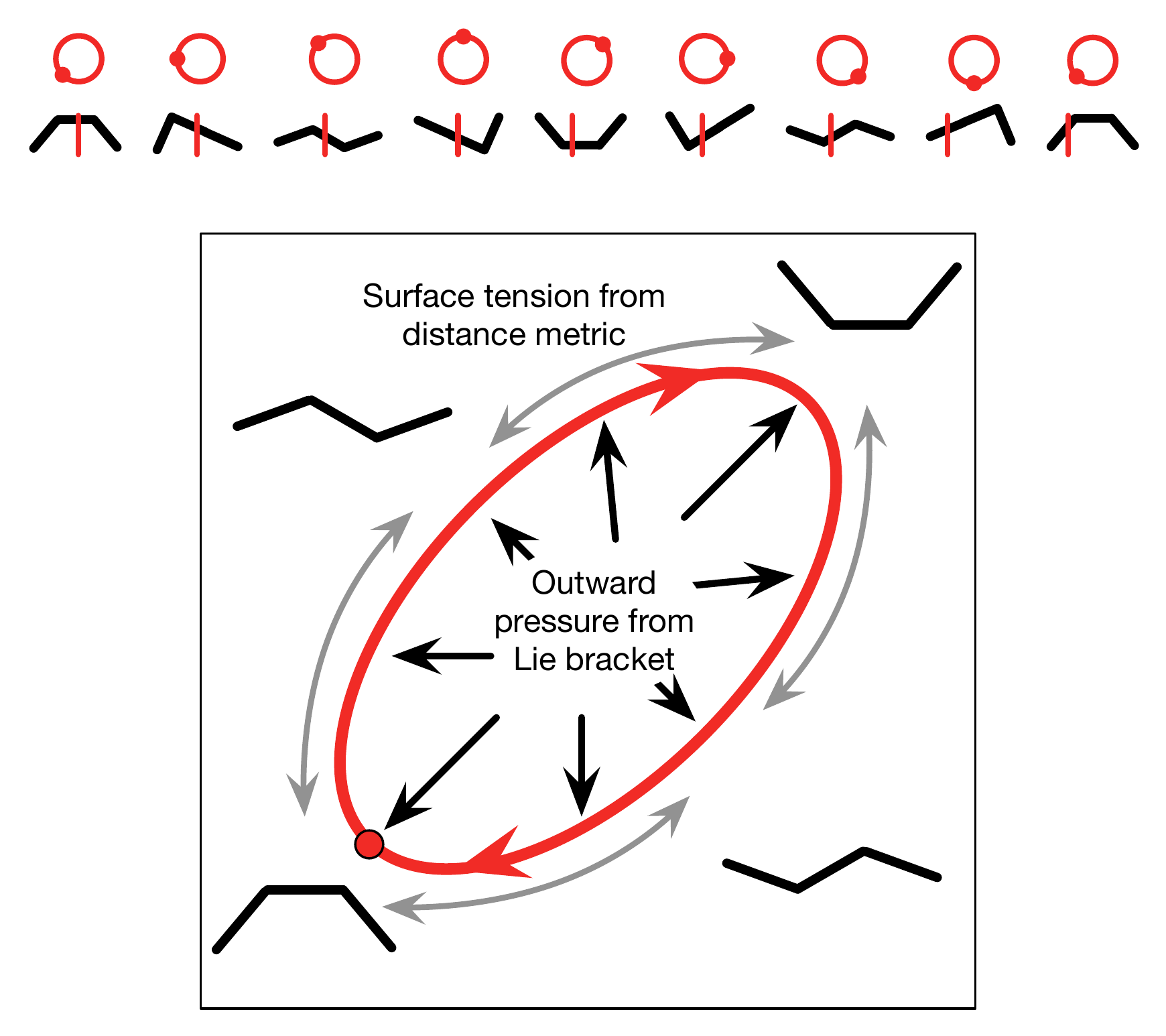}
\caption{Our algorithm maximizes gait efficiency by finding cycles in the space of body shapes that enclose the most \emph{curvature} of the system dynamics (measured via the Lie bracket) while minimizing their cost-to-execute (measured as the metric-weighted lengths of their perimeters). This process is analogous to the process by which air pressure and surface tension combine to produce the shape and size of a soap bubble. Top: The forward progress of a locomoting system as it executes a gait cycle.}
\label{fig:overview}
\end{figure}

An enduring goal~\cite{walsh95,Mukherjee:1993a,Kelly:1995,Radford:1998,Melli:2006,Shammas:2007,Avron:2008} of this community has been to go beyond simply identifying the direction in which small shape oscillations will propel the system, and to use the Lie bracket principle to identify optimal gaits that generate the most displacement, either per gait cycle or per energy expended. For example, a gait that completely encircles a sign-definite region of a height function (by tracing out its zero-contour) would maximize displacement in the corresponding direction; similarly, a contraction of the zero contour that balances the enclosed area against the perimeter length (effort) required to capture it would maximize the efficiency with which this displacement is achieved~\cite{Montgomery:2002vn,Hatton:2011RSS}.

A key obstacle to reaching this goal was the noncommutativity in the equations of motion for most interesting locomoting systems. Because translations and rotations do not commute with each other, summing up Lie brackets over a finite region provides only an approximation of the displacement induced by a gait, and the error in this approximation increases with the amplitude of the oscillation. Historically (e.g.,~\cite{Melli:2006}), this error was regarded as growing too quickly to provide any meaningful information about the optimal gait cycles. In~\cite{Hatton:2009RSS,Hatton:2010ICRA:BVI,Hatton:2011IJRR,Hatton:2015EPJ}, however, we introduced a choice of coordinates that converts much of the system's noncommutativity into nonconservativity, which \emph{is} amenable to finite-scale integration. By working in these new coordinates, we have gained insight into the optimal gaits of swimmers and crawlers in low- and high-viscosity fluids~\cite{Hatton:2013TRO:Swimming,Hatton:2011RSS,Hatton:2011IROS}, and in granular media~\cite{Hatton:2013PRL}. Further, we have used this insight to easily ``hand-draw'' optimized gaits that would otherwise only have been found by a long numerical search of a high-dimensional space of trajectory parameters~\cite{Hatton:2013PRL}.

In this paper, we build on our previous work by encoding the geometric insight described above into a variational gait optimizer. The chief elements of this optimizer are:
\begin{packed_enum}
\ie A gradient ascent/descent component that pushes the cycle to enclose a large sign-definite region of the Lie bracket to maximize the net displacement it generates;
\ie A cost component based on a Riemannian distance metric that limits the growth of the gait cycle; and%
\ie A perimeter-balancing component that evenly-spaces the waypoints in the trajectory, stabilizing the solution and providing an efficiency-optimal parameterization of the resulting motion.
\end{packed_enum}
As illustrated in Fig.~\ref{fig:overview}, these components together form a dynamical system analogous to a soap bubble, whose boundary (here, the gait curve) is pushed outward by internal pressure (here, the Lie bracket) and constrained by surface tension (here, the metric-weighted pathlength). The perimeter-balancing term corresponds to the concentration gradient that evenly distributes soap across the surface of the bubble.

As a demonstration of this approach, we use it to identify optimal gaits for a set of example systems. These systems include Purcell's three-link swimmer~\cite{Purcell:1977} (a standard benchmark for locomotion analysis) and a continuous-curvature extension of Purcell's swimmer~\cite{Hatton:2011RSS} (to emphasize that ``two shape variables'' does not restrict the applicability of this work to systems with ``two joints''). We identify optimal gaits for these systems in both the forward and turning directions, which match those found previously in works such as~\cite{Tam:2007,Hatton:2011IROS}.%

\section{Background}
\label{sec:background}
\subsection{Geometric Locomotion Model}
\label{sec:geo}

\begin{figure}%
\centering
\includegraphics[width=.45\textwidth]{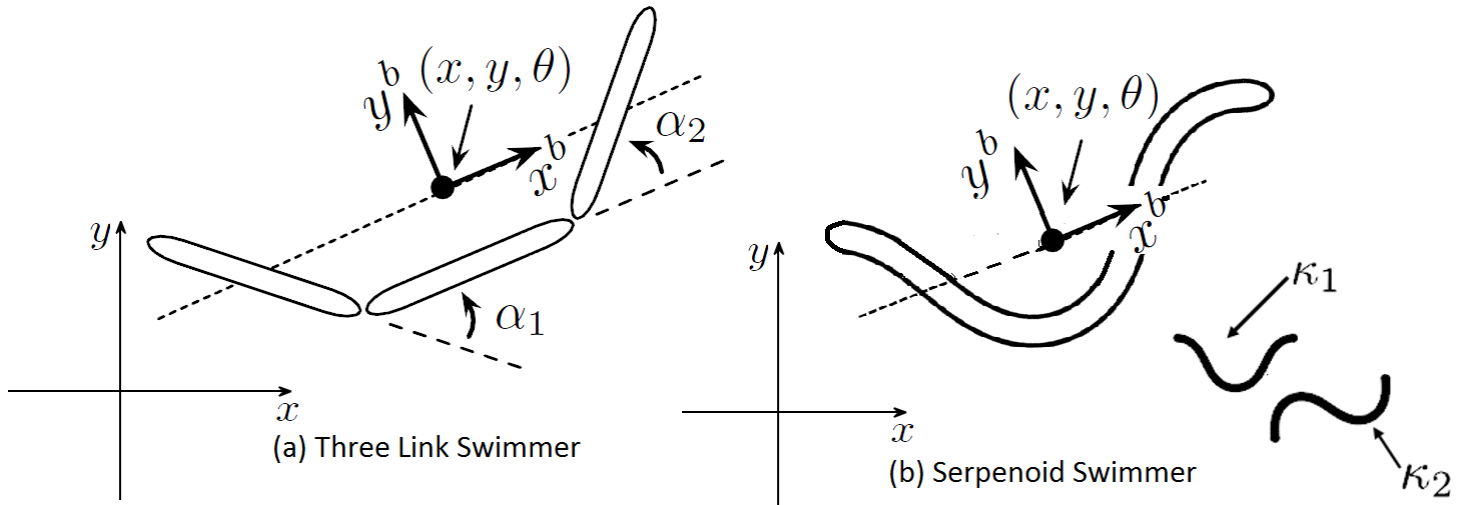}
\caption{Geometry and configuration variables of the two example systems.}
\label{fig:3link}
\end{figure}

When analyzing a locomoting system, it is convenient to separate its configuration space $\bundlespace$ (i.e.\ the space of its generalized coordinates $\bundle$) into a position space $\fiberspace$ and a shape space $\basespace$, such that the position $\fiber\in\fiberspace$ locates the system in the world, and the shape $\base\in\basespace$ gives the relative arrangement of the particles that compose it.\footnote{In the parlance of geometric mechanics, this assigns $\bundlespace$ the structure of a (trivial, principal) \emph{fiber bundle}, with $\fiberspace$ the \emph{fiber space} and $\basespace$ the \emph{base space}.}  For example, the positions of both the three-link system and the serpenoid swimmer in Fig.~\ref{fig:3link} are the locations and orientations of their centroids and mean orientation lines, $\fiber = (x,y,\theta)\in SE(2)$. The shape of the three-link system is parameterized by the two joint angles, $\base = (\alpha_{1},\alpha_{2})$, and the serpenoid's shape is described by a pair of modal amplitudes multiplied by the curvature modes $\kappa_{1}$ and $\kappa_{2}$~\cite{Hatton:2011RSS} that together specify the amplitude and phase of its body wave.

\begin{figure}%
\centering
\includegraphics[scale=.38]{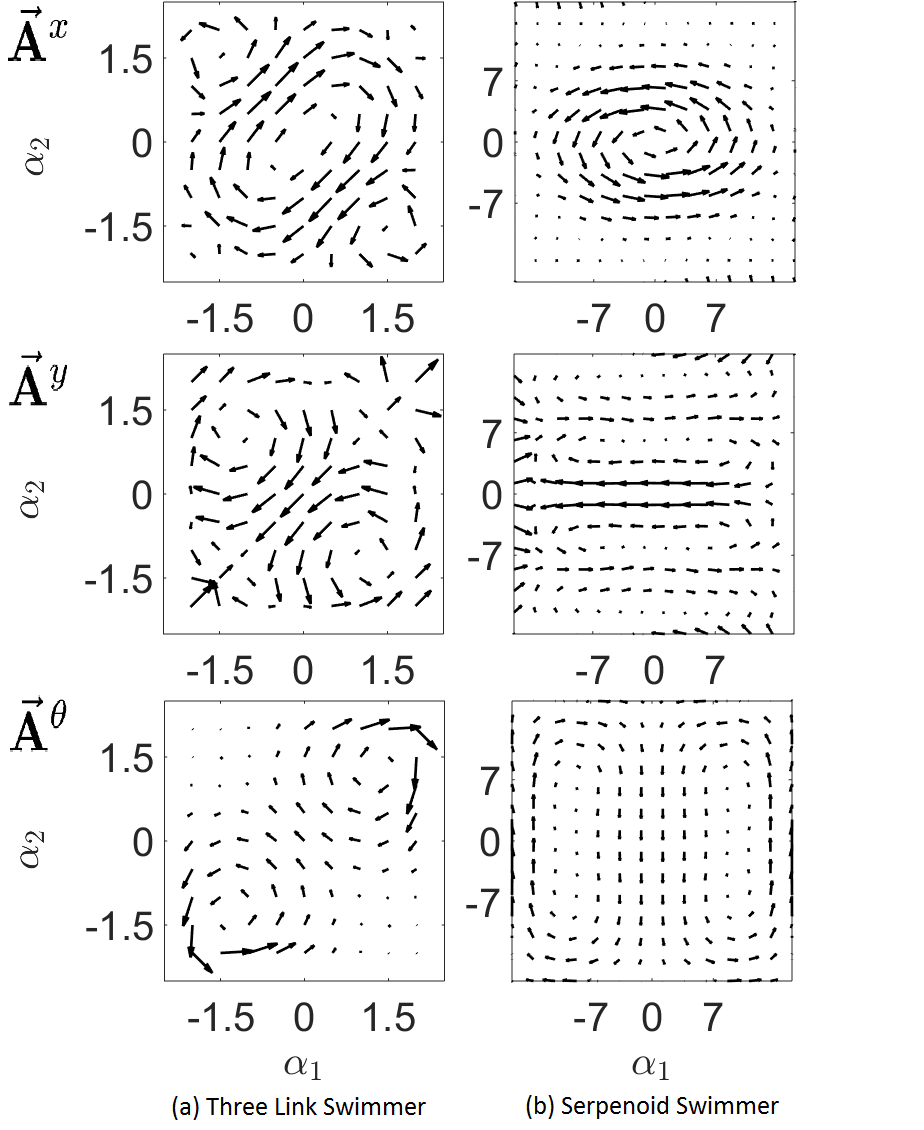}
\caption{The connection vector fields for the three link swimmer is shown in Fig. \ref{fig:connection3link}(a) and that for the serpenoid swimmer are shown in Fig. \ref{fig:connection3link}(b). Note that scale in vector fields has been chosen to emphasize structure, and that scales in different components or systems are not comparable. Also note qualitative similarity (modulo rotation) across systems}
\label{fig:connection3link}
\end{figure}

A useful model for locomotion in \emph{kinematic} regimes where no gliding can occur,\footnote{This kinematic condition has been demonstrated for a wide variety of physical systems, including those whose behavior is dictated by conservation of momentum~\cite{walsh95, Shammas:2007}, non-holonomic constraints such as passive wheels~\cite{Murray:1993, ostrowski98a, Bloch:03,Shammas:2007}, and fluid interactions at the extremes of low~\cite{Avron:2008,Hatton:2010ASME} and high~\cite{Melli:2006,Kanso:2009sw,Hatton:2010ASME} Reynolds numbers.} and  which we employ in this paper, is that at each shape, there exists a linear relationship between changes in the system's shape and changes in its position,
\beq
\bodyvel = - \mixedconn(\base) \basedot,
\label{eq:kinrecon}
\eeq
in which $\bodyvel$ is the body velocity of the system ({i.e.}, $\fiberdot$ expressed in the system's local coordinates), and the \emph{local connection} $\mixedconn$ acts like the Jacobian of a robotic manipulator, mapping from joint velocities to the corresponding body velocity. Each row of $-\mixedconn$ can be regarded as a body-coordinates gradient of one position component with respect to the system shape, as illustrated for the three-link and serpenoid swimmers in Fig.~\ref{fig:connection3link}. %

In a drag-dominated environment, the effort required to change shape can be modeled as the pathlength $s$ of the trajectory through the shape space, weighted by a Riemannian metric $\metric$ as
\begin{equation}\label{eq:distmetric}
d\alnth^{2} = \transpose{d\base} \metric\ d\base.
\end{equation}
As discussed in~\cite{Hatton:2011RSS}, the metric tensor $\metric$ is the same as the map from shape velocity to power dissipated into the surrounding medium,
\beq
P = \metric \basedot,
\eeq
and so can be readily calculated from the first-principles physics of the system

\subsection{Example system dynamics models}

In this paper, we generate the dynamics for our example systems from a resistive force model, in which each element of the body is subject to normal and tangential drag forces proportional to their velocities in those (local) directions. The normal drag coefficient is larger than the tangential component (here, by a factor of $2:1$), corresponding to the general principle that it is harder to move an object in a fluid or on a surface crosswise than it is to move it along its length. We then impose a quasi-static equilibrium condition that the net drag force and moment on the system is zero at all times (treating the system as heavily overdamped, with acceleration forces much smaller than drag forces); because the drag forces are not isometric, the system can use the angle-of-attack of its body surfaces to generate net motion.\footnote{This model is most widely associated with swimmers at low Reynolds numbers (e.g.,~\cite{Tam:2007}),but can also be regarded as an informative general model for systems that experience more lateral drag than longitudinal drag (e.g.~\cite{Hatton:2013PRL}).}

Together, these conditions impose a \emph{Pfaffian constraint}\footnote{A constraint that the allowable velocities are orthogonal to a set of locally-linear constraints, i.e., that they are in the nullspace of a constraint matrix $\omega$.} on the system's generalized velocity,
\beq
\begin{bmatrix} F^{b}_{x} \\ \vspace{2pt} F^{b}_{y} \\ F^{b}_{\theta} \end{bmatrix} = \begin{bmatrix} 0 \\ 0 \\ 0 \end{bmatrix} = \omega(\base) \begin{bmatrix} \bodyvel \\ \dot{\base} \end{bmatrix},\label{eq:forcepfaffian}
\eeq
in which the matrix $\omega$ that maps the velocities to the net forces on the body frame is an expression of the system's internal kinematics and depends only on the shape $\base$. By separating $\omega$ into two sub-blocks, $\omega = [ \omega_{\fiber}^{3\times3}, \omega_{\base}^{3\times n}]$, it is straightforward to rearrange \eqref{eq:forcepfaffian} into
\beq
\bodyvel = -(\omega_{\fiber}^{-1}\omega_{\base})\basedot,\label{eq:invomega}
\eeq
revealing the local connection as $\mixedconn = \omega_{\fiber}^{-1}\omega_{\base}$. For a more detailed treatment this process, we refer the reader to~\cite{Hatton:2013TRO:Swimming}.  Once $\mixedconn$ has been found, it can be used to calculate the Riemannian metric $\metric$ as
\beq
\metric = \int_{\text{body}} J^{T}(\ell)\ c\ J(\ell) \ d\ell,
\eeq
where $J(\ell)$ is the Jacobian from shape velocity to the local velocity of each section of the body (which incorporates both $\mixedconn$ and the system's internal kinematics), and $c$ is the matrix of drag coefficients.

We apply this physics models to two example geometries: Purcell's three-link swimmer~\cite{Purcell:1977}, and a \emph{serpenoid swimmer}~\cite{hirose93}. The three-link swimmer is a useful and widely-adopted~\cite{Tam:2007,Avron:2008,Melli:2006,Hatton:2013TRO:Swimming,Bettiol:2015,Giraldi:2013aa,DeSimone:2012aa} minimal example for locomotion, because the two degrees of freedom can be easily visualized. The serpenoid swimmer, whose shape is defined by the amplitude of sine and cosine curvature modes~\cite{Hatton:2011RSS}, provides an example of a two-DoF system that has been shown to closely model how snakes and other animals use undulatory locomotion to move through the world~\cite{hirose93}.

\subsection{Finding optimal gaits}\label{sec:gaitoptimization}

Optimal gait design has a long history of research in the physics, mathematics, and engineering communities, as part of the broader field of optimal control~\cite{Bryson:1979,Ostrowski:2000vl}. For systems of the classes we consider here, notable contributions include those of Purcell, who introduced the three-link swimmer as a minimal template for understanding locomotion, a series of works~\cite{Tam:2007,DeSimone:2012aa,Giraldi:2013aa,Bettiol:2015} aimed at numerically optimizing the stroke pattern, and the observation in~\cite{Becker:2003} that the optimal pacing for the gait keeps the power dissipation constant over the cycle.

\begin{figure}%
\centering
\includegraphics[scale=0.38]{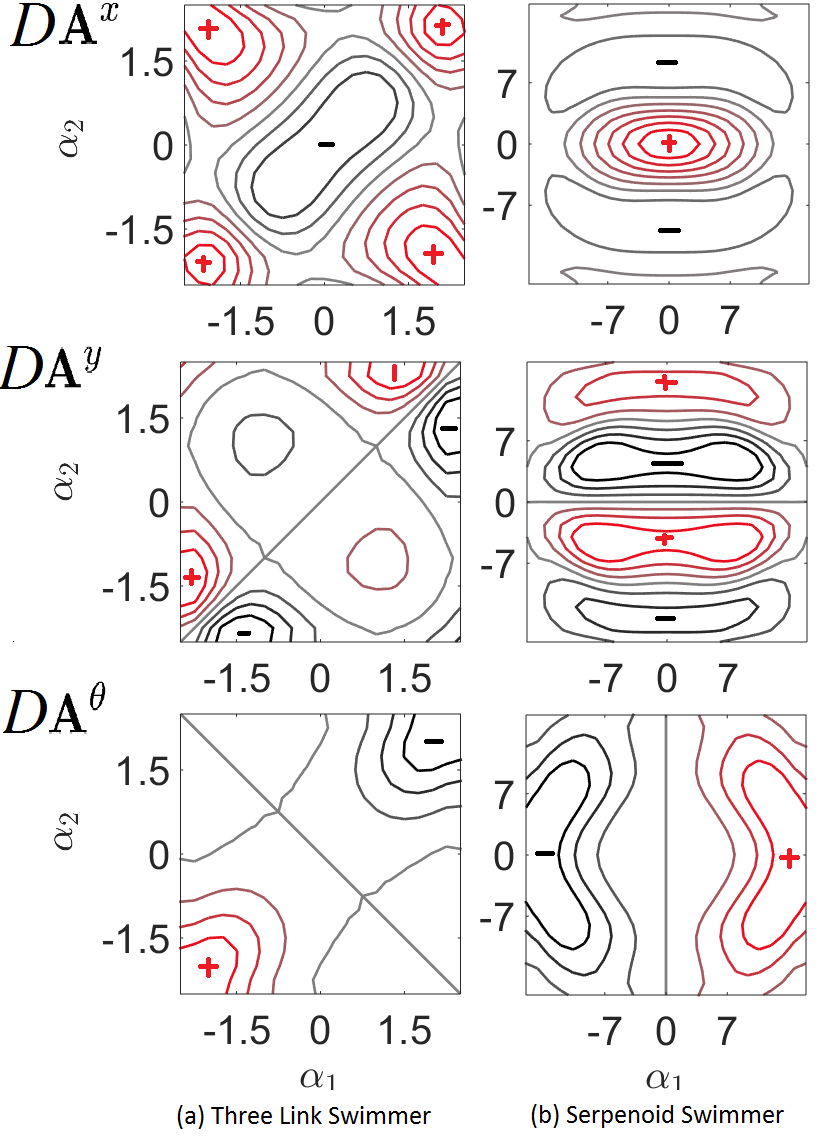}
\caption{Height Functions (plots of the Lie bracket $D\mixedconn$) for the three link swimmer are shown in Fig.\ref{fig:heightfun}(a) and that for the serpenoid swimmer are shown in Fig.\ref{fig:heightfun}(b). Note the qualitative similarity (modulo a rotation) between the two sets of height functions.}
\label{fig:heightfun}
\end{figure}

Several efforts  in the geometric mechanics community (including our own), have aimed to use the structure of the systems' \emph{Lie brackets} (a measure of how ``non-canceling" the system dynamics are over infinitesimal cyclic inputs) to understand the structure of the optimal solutions to the system equations of motion~\cite{walsh95,ostrowski98a,Melli:2006,Shammas:2007,Avron:2008,Hatton:2013TRO:Swimming,Hatton:2011RSS}. The core principle in these works is that because the net displacement $\gaitdisp$ over a gait cycle~$\gait$ is the line integral of~\eqref{eq:kinrecon} along $\gait$, the displacement can be approximated\footnote{This approximation (a generalized form of Stokes' theorem) is a truncation of the Baker-Campbell-Hausdorf series for path-ordered exponentiation on a noncommutative group, and closely related to the Magnus expansion~\cite{Radford:1998,Magnus:1954vl}. For a discussion of the accuracy of this approximation and its derivation, see~\cite{Hatton:2015EPJ}. In presenting this approximation, we also elide some details of exponential coordinates on Lie groups, which are also discussed in~\cite{Hatton:2015EPJ}.} by an area integral of the curvature $D\mixedconn$ of the local connection (its total Lie bracket~\cite{Hatton:2015EPJ})   over a surface $\gait_{a}$ bounded by the cycle:
\begin{align}
\gaitdisp &= \ointctrclockwise_{\gait} - \fiber\mixedconn(\base) \label{eq:gaitpathintegral}\\ &\approx \iint_{\gait_{a}} \underbrace{-\extd\mixedconn + \big{[}\mixedconn_{1},\mixedconn_{2}\big{]}}_{\text{$-D\mixedconn$ (total Lie bracket)}}, \label{eq:lie}
\end{align}
in which $\extd\mixedconn$ is the exterior derivative of the local connection (its generalized row-wise curl) and the local Lie bracket term evaluates (on $SE(2)$) as
\beq
\big{[}\mixedconn_{1},\mixedconn_{2}\big{]} = \begin{bmatrix} \mixedconn^{y}_{1}\mixedconn^{\theta}_{2}- \mixedconn^{y}_{2}\mixedconn^{\theta}_{1} \\ \mixedconn^{x}_{2}\mixedconn^{\theta}_{1}- \mixedconn^{x}_{1}\mixedconn^{\theta}_{2} \\ 0 \end{bmatrix}.
\eeq

Plotting these curvature terms as height functions (as in Fig.~\ref{fig:heightfun}) reveals the attractors that influence the optimal gait cycles: Gaits that produce net displacement in a given $(x,y,\theta)$ direction are located in strongly sign-definite regions of the corresponding $D\mixedconn$ height function. For example, $x$-translation gaits encircle the center of the shape space for both the three-link and serpenoid systems, whereas $y$-translations or $\theta$-rotations are produced by cycles in the corners or edges of the shape space.

In~\cite{Hatton:2010ICRA:BVI,Hatton:2011IJRR,Hatton:2015EPJ}, we identified coordinate choices that make the approximation in~\eqref{eq:lie} accurate for large-amplitude gaits, which allowed us to make two extensions to the area principles: First, that maximum-displacement gaits (the ``longest strides'' that the systems can take) follow the zero-contours of the height functions, completely enclosing a sign definite region. Second, that the most efficient gaits (the systems' ``comfortable strides'') are contractions/straightenings of these zero-contours; by giving up the low-yield regions at the edges of the sign-definite regions (or crossing slightly into opposite-sign regions) the system travels a shorter distance in each cycle, but the consequently shorter perimeter length means the system can repeat the cycles more quickly at a given level of power consumption.

\section{Optimizer}

Our aim with this paper is to encode the geometric principles described in~\S\ref{sec:gaitoptimization} into a gait optimization algorithm.

We start from the basic variational principle that functions reach their extrema when their derivatives go to zero. Given a gait parameterization $\varparam$, maximum-displacement cycles therefore satisfy the condition that the gradient of net displacement with respect to the parameters is zero,
\beq \label{eq:maxdispcond}
\nabla_{\varparam}\gaitdisp = \mathbf{0}
\eeq
and maximum-efficiency cycles (normalizing displacement by the pathlength-effort required to generate it) likewise satisfy the condition that the gradient of this efficiency ratio is zero,
\beq \label{eq:maxeffcond}
\nabla_{\varparam}\frac{\gaitdisp}{\alnth} = \frac{1}{\alnth} \nabla_{\varparam}\gaitdisp -\frac{\gaitdisp}{\alnth^{2}}\nabla_{\varparam}\alnth= \mathbf{0}.
\eeq
For suitable seed values $\varparam_{0}$, solutions to~\eqref{eq:maxdispcond} and~\eqref{eq:maxeffcond} can therefore be reached by finding the respective equilibria of the dynamical systems
\beq \label{eq:gaitequilibria}
\dot{\varparam} = \nabla_{\varparam}\gaitdisp \ \ \text{and} \ \ \dot{\varparam} = \nabla_{\varparam}{\frac{\gaitdisp}{\alnth}}.
\eeq

As discussed in the following subsections, we then incorporate the geometric principles in~\eqref{eq:distmetric} and~\eqref{eq:lie} to evaluate the $\nabla_{\varparam}\gaitdisp$ and $\nabla_{\varparam}\alnth$ terms in~\eqref{eq:maxeffcond} as functions of the system's Lie bracket $D\mixedconn$ and Riemannian metric $\metric$. We also use the metric $\metric$ to introduce a third term,~$\nabla_{\varparam}\intstress$, to~\eqref{eq:gaitequilibria}. This term acts orthogonally to the optimization process, but serves to enforce an even sampling of points around the perimeter of the gait. This even distribution promotes stability of the optimizer and additionally corresponds to optimal pacing (internal timing) of the gait, as discussed in~\cite{Becker:2003,Tam:2007,Hatton:2011RSS}.

Collecting all three terms and factoring out a coefficient of~$\frac{1}{\alnth}$ from~\eqref{eq:maxeffcond}, our algorithm thus finds the maximum-efficiency gait as the equilibrium of
\beq \label{eq:fullequation}
\dot{\varparam} = \nabla_{\varparam}\gaitdisp -\frac{\gaitdisp}{\alnth}\nabla_{\varparam}\alnth + \nabla_{\varparam}\intstress.
\eeq
As illustrated in Fig.~\ref{fig:overview}, this differential equation is directly analogous to the equations governing the shape of a soap bubble: $\nabla_{\varparam}\gaitdisp$ takes the Lie bracket as an ``internal pressure'' seeking to expand the gait cycle, $\nabla_{\varparam}\alnth$ is the ``surface tension'' that constrains the growth of the bubble, and $\nabla_{\varparam}\intstress$ is the ``concentration gradient'' that spreads the soap over the bubble's surface.

In the following subsections, we explore each of the terms in~\eqref{eq:fullequation}, discussing both their fundamental geometric definitions and how they can be implemented in a direct-transcription solver that parameterizes the gait as a sequence of waypoints.

\subsection{Internal Pressure from the Lie Bracket}

The first term in~\eqref{eq:fullequation},~$\nabla_{\varparam}\gaitdisp$, guides the gait towards maximum-displacement cycles. By sustituting the approximation from~\eqref{eq:lie} into this expression as
\beq
\nabla_{\varparam}\gaitdisp \approx \nabla_{\varparam}\iint_{\gait_{a}}{(-D\mixedconn)},
\eeq
and noting that variations in $\varparam$ affect $\gait$ but not $D\mixedconn$, we can convert~$\nabla_{\varparam}\gaitdisp$ into gradient of an area integral with respect to variations in its boundary. We can then invoke a powerful geometric principle,\footnote{The general form of the \emph{Leibniz integral rule}~\cite{Flanders:1973aa}.} which states

\emph{
The gradient of an integral with respect to variations of its boundary is equal to the gradient of the boundary with respect to these variations, multiplied by the integrand evaluated along the boundary.
}

Formally, this multiplication is the \emph{interior product}\footnote{Not the inner product, see~\cite{Flanders:1973aa} for more details.} of the boundary gradient with the integrand,
\beq \label{eq:interiorproduct}
\nabla_{\varparam}\iint_{\gait_{a}}{(-D\mixedconn)} = \ointctrclockwise_{\gait} (\nabla_{\varparam}\gait) \intprod (-D\mixedconn),
\eeq
which contracts $D\mixedconn$ (a \emph{differential two-form}~\cite{Kobayashi:1963}) along $\nabla_{\varparam}\gait$ to produce a differential one-form that can be integrated over $\gait$. This formalism will become important in future work when we explore these principles on systems with more than two shape variables; for the two-dimensional shape spaces we consider in this paper, the interior product reduces to a simple multiplication between the outward component of $\nabla_{\varparam}\gait$ and the scalar magnitude of the Lie bracket,
\beq \label{eq:boundarygradientsimple}
\nabla_{\varparam}\iint_{\gait_{a}}{(-D\mixedconn)} = \ointctrclockwise_{\gait} (\nabla_{\varparam_\perp}\gait) (-D\mixedconn).
\eeq

\begin{figure*}%
\centering
\includegraphics[width=.92\textwidth]{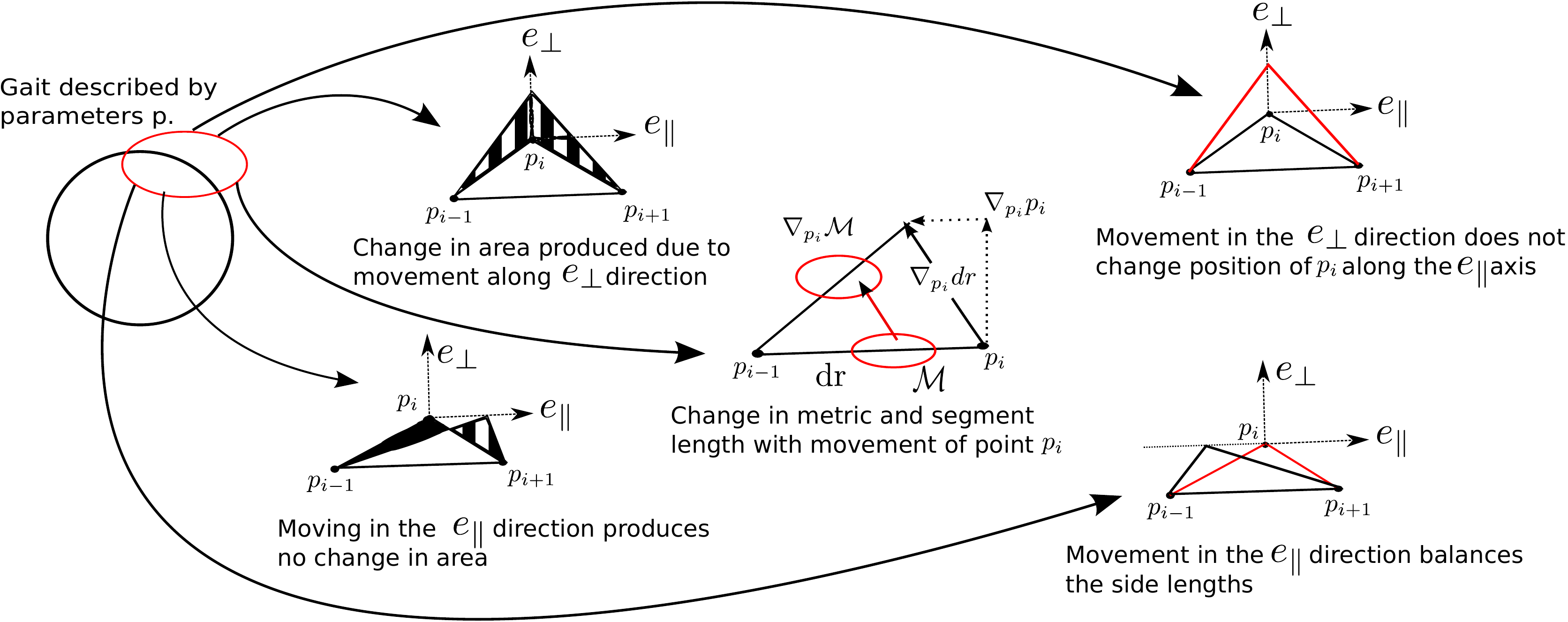}
\caption{This figure illustrates the changes in area caused by moving in the two coordinate directions in the local frame. Moving in the tangential direction $\localbasis_{\parallel}$ produces no change in area, as the area of the triangle given by half the product of base length and height remains the same.}
\label{fig:derivation}
\end{figure*}

\emph{Implementation of the internal pressure:}
Each waypoint $\varparam_{i}$ forms a triangle with its neighboring points, whose base defines a local tangent direction $\localbasis_{\parallel}$ as
\beq
\varparam_{i+1} - \varparam_{i-1} = \ell\, \localbasis_{\parallel}
\eeq
and a local normal direction $\localbasis_{\perp}$ orthogonal to $\localbasis_{\parallel}$. As illustrated in Fig.~\ref{fig:derivation}, the gradient of the enclosed area with respect to variations in the position of $\varparam_{i}$ in the~$\localbasis_{\parallel}$ and~$\localbasis_{\perp}$ directions is the change in triangle's area as $\varparam_{i}$ moves. Because the triangle's area is always one half base times height (regardless of its pitch or the ratio of its sidelengths), this gradient evaluates to
\beq
\nabla_{\varparam_{i}}\gait_{a} =  \begin{bmatrix} \localbasis_{\parallel} & \localbasis_{\perp} \end{bmatrix}\begin{bmatrix} 0 \\ \ell/2 \end{bmatrix}.
\eeq
Note that this term matches the right-hand side of~\eqref{eq:boundarygradientsimple}, with only normal motions of the boundary affecting the enclosed area.

\subsection{Surface Tension from the Distance Metric}
The second term in~\eqref{eq:fullequation} takes $\nabla_{\varparam}\alnth$ as a measure of how variations in the gait affect the cost of executing it, and scales this by a factor of $\frac{\gaitdisp}{\alnth}$ to compare how the return on this investment compares to the efficiency of the gait in its present state. The gradient component of this term can be related to the system's Riemannian metric by first incorporating the arclength calculation from~\eqref{eq:distmetric},
\beq
\nabla_{\varparam}\alnth = \nabla_{\varparam}\oint_{\gait} (\overbrace{\transpose{d\base} \metric\ d\base}^{d\alnth^{2}})^{\frac{1}{2}}
\eeq
and then applying standard calculus operations\footnote{Namely: Differentiation under the integral sign, chain rule, product rule, and then exploiting the symmetry of $\metric$ to consolidate terms.} to arrive at
\beq \label{eq:surfacetensionfull}
\nabla_{\varparam}\alnth = \frac{1}{2\alnth}\oint_{\gait} \Big( 2\transpose{(\nabla_{\varparam}d\base)} \metric\ d\base + \transpose{d\base} (\nabla_{\varparam}\metric)\ d\base \Big),
\eeq
in which the two parts of the integrand respectively measure how changes in the relative positions of the boundary elements and changes in the metric at the underlying points affect the pathlength, and hence the cost of motion.

The $\frac{\gaitdisp}{\alnth}$ factor (which normalizes the scales of $\nabla_{\varparam}\gaitdisp$ and $\nabla_{\varparam}\alnth$) can be calculated directly from~\eqref{eq:gaitpathintegral} and the integral of~\eqref{eq:distmetric}. Although the calculation of $\gaitdisp$ could in theory make use of the area approximation in~\eqref{eq:lie}, this would be inefficient and impractical: integration of surfaces with arbitrarily complex boundaries requires significantly more computational resources than are needed for line integration around the boundary. Using the true line integral also improves the accuracy of the solution; by continuously recalibrating to the true net displacement, the algorithm avoids compounding any errors introduced by the approximation in~\eqref{eq:lie}.

\emph{Implementation of the surface tension:}
Each waypoint $\varparam_i$ is at the head of a vector extending from $\varparam_{i-1}$, such that the $d\base$ vector and its gradient in~\eqref{eq:surfacetensionfull} can be taken as
\beq
d\base = \varparam_{i} - \varparam_{i-1} \ \ \ \text{and} \ \ \ \nabla_{\varparam_{i}}d\base = \begin{bmatrix} 1 \\ 1 \end{bmatrix}.
\eeq
For computational simplicity, we evaluate the metric at the center of each segment. This point moves at the mean speed of the segment endpoints, and its gradient with respect to changes in $\varparam_{i}$ is thus one half of its gradient over the underlying space,
\beq
\nabla_{\varparam_{i}} \metric = \frac{1}{2} \nabla \metric.
\eeq
These gradient relationships are illustrated in Fig.~\ref{fig:derivation}(b), with the metric represented by its Tissot indicatrix ellipse~\cite{Hatton:2011RSS}.

\subsection{Concentration Gradient from Parameterization}

The third term in~\eqref{eq:fullequation}, $\nabla_{\varparam}\intstress$, is included to ensure an even sampling of points around the gait during numerical optimization. This term does not appear in~\eqref{eq:maxdispcond} or~\eqref{eq:maxeffcond} because it is orthogonal to those optimization criteria: varying a continuous curve along itself does not change its length or the region it encloses. When the curve is represented with a discrete basis, however, that basis induces a length coordinate~$\altlnth$ along the gait that is not necessarily compatible with the natural coordinate $\alnth$. This coordinate may lead to undersampling of some portions of the curve or pathological behavior, such as crossovers in which the gait spontaneously forms internal loops that do not correspond to features in $D\mixedconn$.

To reduce the chance of such instabilities occurring, we first identify an internal strain energy~$\intstress$ that measures how much an even spacing according to~$\altlnth$ differs from even spacing according to~$\alnth$. This strain energy,
\beq \label{eq:intstress}
\sigma = \oint_{\gait} \left(\frac{\alnth}{\altlnth} - \nabla_{\altlnth} \alnth\right)^{2},
\eeq
compares the average and local relative rates at which the two coordinates are changing around the curve; if these rates are not matched, then the two measures of length are imbalanced. Including the gradient term of this energy with respect to variations in the parameters, $\nabla_{\varparam}\intstress$, then directs the optimizer to use its degrees of freedom that are (instantaneously) in the kernel of the primary optimization criteria to reduce the imbalance. The expansion of this term depends on the choice of parameters for the gait (because they affect the construction of~$\altlnth$), so we do not provide a general formula for it. In the case of direct transcription, however, it works out to be very straightforward.

\emph{Implementation of the concentration gradient:}
In a direct-transcription parameterization, the strain energy at a given waypoint corresponds to the difference in the tangential distance from that waypoint to each of its neighbors,
\beq
\intstress_{i} = \big((\varparam_{i+1}-\varparam_{i})_{\parallel} - (\varparam_{i}-\varparam_{i-1})_{\parallel}\big)^{2},
\eeq
whose gradient with respect to the position of $\varparam_{i}$ is proportional to the sum of the tangential displacements of the neighboring points relative to $\varparam_{i}$,
\beq
\nabla_{\varparam_{i}}\intstress \propto  \begin{bmatrix} \localbasis_{\parallel} & \localbasis_{\perp} \end{bmatrix}\begin{bmatrix} (\varparam_{i+1})_{\parallel} + (\varparam_{i-1})_{\parallel} \\ 0 \end{bmatrix}.
\eeq

\section{Demonstration}

We implemented the optimizer described in the previous section in Matlab, both directly solving the differential equation in~\eqref{eq:fullequation} using \verb#ode45# and by providing~\eqref{eq:fullequation} as the gradient for the standard \verb#fmincon# optimizer using the interior-point algorithm. As expected, both implementations converged on the same solutions, with the \verb#fmincon# implementation completing more quickly (on the order of minutes on a modern desktop computer for a gait with $100$ waypoints) due to its ability to take larger steps through the parameter space.%

\begin{figure*}
\centering
\includegraphics[width=.9\textwidth]{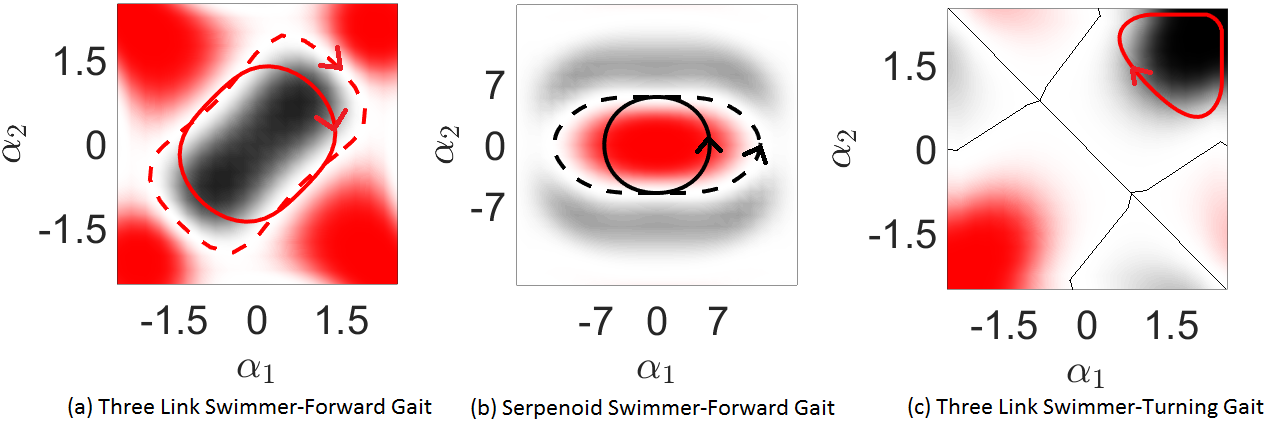}
\caption{Optimal gaits found by our algorithm. Maximum-displacement gaits, which follow the zero-contours of the corresponding height function, are indicated by dashed lines. Maximum-efficiency gaits, which are contractions of the zero-contour, are indicated by solid lines. Orientation on the curves is chosen so that the displacement from each gait is positive.}
\label{fig:Optimized gait}
\end{figure*}

We ran the optimizer on the three link system twice, first to find the gait that maximizes the displacement in the $x$-direction over a single cycle and then to find the maximum-efficiency cycle. The dashed line in Fig.~\ref{fig:Optimized gait}(a) shows the gait that optimizes the maximum displacement over the cycle for the three-link swimmer. As expected, this gait follows the zero-contour of the height function. The solid line Fig.~\ref{fig:Optimized gait}(a) shows the maximum efficiency gait for this system. Because the maximum-efficiency optimizer places a cost on pathlength, this curve gives up the low-yield regions at the ends of the cycle and crosses slightly outside the zero contour.

Fig.~\ref{fig:Optimized gait}(b) illustrates essentially-similar behavior for the serpenoid: the dashed maximum-displacement gait traces the zero contour, and the solid maximum-efficiency gait captures a more compact area within the sign-definite region.
Fig.~\ref{fig:Optimized gait}(c) shows the gait that produces the most cost-effective rotational motion of the three-link swimmer (the most sign-definite area it could capture on the $\theta$ height function without over-stretching itself). The results the optimizer yields for the three link swimmer match those obtained  for maximum efficiency in~\cite{Hatton:2011IROS} and, along the $x$ direction, those in~\cite{Tam:2007}.

\section{Conclusions}

In this paper, we have incorporated recent geometric insights about optimal gaits for kinematic systems into a variational trajectory optimizer. The dynamics of this optimizer are analogous to those of a soap bubble, with the Lie bracket providing an ``inflating pressure'' to the trajectory and the Riemannian metric on the shape space contributing ``surface tension'' that halts growth of the cycle in the face of diminishing returns. Together, these elements drive the gait cycle to a ``comfortable stride'' that converts shape change effort into net displacement with optimal efficiency.

We demonstrated this optimizer in operation on a number of test cases, including Purcell's three-link swimmer (a standard minimal template for locomotion modeling) and a serpenoid swimmer (a model widely used in studies of animals and snake robots). The optimizer correctly found the maximum-displacement and maximum efficiency cycles for these systems, as corroborated by previous exhaustive optimizations of the gait cycles such as~\cite{Tam:2007}.

A line of future work that we are pursuing is to extend this optimizer to systems whose dynamics are dominated by inertial effects, including those subject to conservation of momentum whose momentum is directed by nonholonomic constraints. The cost of a gait for these types of systems depends on the acceleration through the configuration shape during the gait rather than on the energy dissipated into the environment. Therefore we expect that a cost associated with the gait's curvature through the shape space (as well as its length) will play an important role in finding the optimal gait for these systems.

\bibliographystyle{IEEEtran}

\begin{thebibliography}{}
\providecommand{\url}[1]{#1}
\csname url@samestyle\endcsname
\providecommand{\newblock}{\relax}
\providecommand{\bibinfo}[2]{#2}
\providecommand{\BIBentrySTDinterwordspacing}{\spaceskip=0pt\relax}
\providecommand{\BIBentryALTinterwordstretchfactor}{4}
\providecommand{\BIBentryALTinterwordspacing}{\spaceskip=\fontdimen2\font plus
\BIBentryALTinterwordstretchfactor\fontdimen3\font minus
  \fontdimen4\font\relax}
\providecommand{\BIBforeignlanguage}[2]{{%
\expandafter\ifx\csname l@#1\endcsname\relax
\typeout{** WARNING: IEEEtran.bst: No hyphenation pattern has been}%
\typeout{** loaded for the language `#1'. Using the pattern for}%
\typeout{** the default language instead.}%
\else
\language=\csname l@#1\endcsname
\fi
#2}}
\providecommand{\BIBdecl}{\relax}
\BIBdecl

\end{thebibliography}


\begin{thebibliography}{10}
\providecommand{\url}[1]{#1}
\csname url@samestyle\endcsname
\providecommand{\newblock}{\relax}
\providecommand{\bibinfo}[2]{#2}
\providecommand{\BIBentrySTDinterwordspacing}{\spaceskip=0pt\relax}
\providecommand{\BIBentryALTinterwordstretchfactor}{4}
\providecommand{\BIBentryALTinterwordspacing}{\spaceskip=\fontdimen2\font plus
\BIBentryALTinterwordstretchfactor\fontdimen3\font minus
  \fontdimen4\font\relax}
\providecommand{\BIBforeignlanguage}[2]{{%
\expandafter\ifx\csname l@#1\endcsname\relax
\typeout{** WARNING: IEEEtran.bst: No hyphenation pattern has been}%
\typeout{** loaded for the language `#1'. Using the pattern for}%
\typeout{** the default language instead.}%
\else
\language=\csname l@#1\endcsname
\fi
#2}}
\providecommand{\BIBdecl}{\relax}
\BIBdecl

\bibitem{Murray:1993}
R.~M. Murray and S.~S. Sastry, ``Nonholonomic motion planning: Steering using
  sinusoids,'' \emph{IEEE Transactions on Automatic Control}, vol.~38, no.~5,
  pp. 700--716, Jan 1993.

\bibitem{Morgansen:2007}
\BIBentryALTinterwordspacing
K.~A. Morgansen, B.~I. Triplett, and D.~J. Klein, ``Geometric methods for
  modeling and control of free-swimming fin-actuated underwater vehicles,''
  \emph{IEEE Transactions on Robotics}, vol.~23, no.~6, pp. 1184--1199, Jan
  2007. [Online]. Available:
  \url{http://ieeexplore.ieee.org/xpls/abs_all.jsp?arnumber=4399955}
\BIBentrySTDinterwordspacing

\bibitem{walsh95}
G.~C. Walsh and S.~Sastry, ``On reorienting linked rigid bodies using internal
  motions,'' \emph{Robotics and Automation, IEEE Transactions on}, vol.~11,
  no.~1, pp. 139--146, January 1995.

\bibitem{ostrowski98a}
J.~P. Ostrowski and J.~Burdick, ``The mechanics and control of undulatory
  locomotion,'' \emph{International Journal of Robotics Research}, vol.~17,
  no.~7, pp. 683--701, July 1998.

\bibitem{Shammas:2007}
E.~A. Shammas, H.~Choset, and A.~A. Rizzi, ``Geometric motion planning analysis
  for two classes of underactuated mechanical systems,'' \emph{Int. J. of
  Robotics Research}, vol.~26, no.~10, pp. 1043--1073, 2007.

\bibitem{Avron:2008}
J.~E. Avron and O.~Raz, ``A geometric theory of swimming: {P}urcell's swimmer
  and its symmetrized cousin,'' \emph{New Journal of Physics}, vol.~9, no. 437,
  2008.

\bibitem{Hatton:2013TRO:Swimming}
R.~L. Hatton and H.~Choset, ``Geometric swimming at low and high reynolds
  numbers,'' \emph{IEEE Transactions on Robotics}, vol.~29, no.~3, pp.
  615--624, June 2013.

\bibitem{Melli:2006}
J.~B. Melli, C.~W. Rowley, and D.~S. Rufat, ``Motion planning for an
  articulated body in a perfect planar fluid,'' \emph{SIAM Journal of Applied
  Dynamical Systems}, vol.~5, no.~4, pp. 650--669, November 2006.

\bibitem{Kanso:2009sw}
E.~Kanso, ``Swimming due to transverse shape deformations,'' \emph{Journal of
  Fluid Mechanics}, vol. 631, pp. 127--148, 2009.

\bibitem{Hatton:2013PRL}
R.~L. Hatton, H.~Choset, Y.~Ding, and D.~I. Goldman, ``Geometric visualization
  of self-propulsion in a complex medium,'' \emph{Physical Review Letters},
  vol. 110, p. 078101, February 2013.

\bibitem{Mukherjee:1993a}
R.~Mukherjee and D.~P. Anderson, ``A surface integral approach to the motion
  planning of nonholonomic systems,'' in \emph{American Control Conference,
  1993}, 1993, pp. 1816 --1823.

\bibitem{Kelly:1995}
S.~D.~K. Kelly and R.~M. Murray, ``Geometric phases and robotic locomotion,''
  \emph{J. Robotic Systems}, vol.~12, no.~6, pp. 417--431, Jan 1995.

\bibitem{Radford:1998}
J.~E. Radford and J.~W. Burdick, ``Local motion planning for nonholonomic
  control systems evolving on principal bundles,'' in \emph{{Proceedings of the
  International Symposium on Mathematical Theory of Networks and Systems}},
  Padova, Italy, 1998.

\bibitem{Montgomery:2002vn}
R.~Montgomery, \emph{A Tour of Subriemannian Geometries, Their Geodesics, and
  Applications}.\hskip 1em plus 0.5em minus 0.4em\relax American Mathematical
  Society, 2002.

\bibitem{Hatton:2011RSS}
R.~L. Hatton and H.~Choset, ``Kinematic cartography for locomotion at low
  {R}eynolds numbers,'' in \emph{{Proceedings of Robotics: Science and Systems
  VII}}, Los Angeles, CA USA, June 2011.

\bibitem{Hatton:2009RSS}
------, ``Approximating displacement with the body velocity integral,'' in
  \emph{{Proceedings of Robotics: Science and Systems V}}, Seattle, WA USA,
  June 2009.

\bibitem{Hatton:2010ICRA:BVI}
------, ``Optimizing coordinate choice for locomoting systems,'' in
  \emph{{Proceedings of the IEEE International Conference on Robotics and
  Automation}}, Anchorage, AK USA, May 2010, pp. 4493--4498.

\bibitem{Hatton:2011IJRR}
------, ``Geometric motion planning: The local connection, {S}tokes' theorem,
  and the importance of coordinate choice,'' \emph{International Journal of
  Robotics Research}, vol.~30, no.~8, pp. 988--1014, July 2011.

\bibitem{Hatton:2015EPJ}
------, ``Nonconservativity and noncommutativity in locomotion,''
  \emph{European Physical Journal Special Topics: Dynamics of Animal Systems},
  vol. 224, no. 17--18, pp. 3141--3174, 2015.

\bibitem{Hatton:2011IROS}
R.~L. Hatton, L.~J. Burton, A.~E. Hosoi, and H.~Choset, ``Geometric
  maneuverability, with applications to low {R}eynolds number swimming,'' in
  \emph{{Proceedings of the IEEE/RSJ International Conference on Intelligent
  Robots and Systems}}, San Francisco, CA USA, September 2011.

\bibitem{Purcell:1977}
E.~M. Purcell, ``Life at low {R}eynolds numbers,'' \emph{American Journal of
  Physics}, vol.~45, no.~1, pp. 3--11, January 1977.

\bibitem{Tam:2007}
D.~Tam and A.~E. Hosoi, ``Optimal stroke patterns for {P}urcell's three-link
  swimmer,'' \emph{Phys. Review Letters}, vol.~98, no.~6, p. 068105, 2007.

\bibitem{Bloch:03}
A.~M. Bloch \emph{et~al.}, \emph{Nonholonomic Mechanics and Control}.\hskip 1em
  plus 0.5em minus 0.4em\relax Springer, 2003.

\bibitem{Hatton:2010ASME}
R.~L. Hatton and H.~Choset, ``Connection vector fields and optimized
  coordinates for swimming systems at low and high {R}eynolds numbers,'' in
  \emph{{Proceedings of the ASME Dynamic Systems and Controls Conference
  (DSCC)}}, Cambridge, Massachusetts, USA, Sep 2010.

\bibitem{hirose93}
S.~Hirose, \emph{Biologically Inspired Robots (Snake-like Locomotor and
  Manipulator)}.\hskip 1em plus 0.5em minus 0.4em\relax Oxford University
  Press, 1993.

\bibitem{Bettiol:2015}
\BIBentryALTinterwordspacing
P.~Bettiol, B.~Bonnard, L.~Giraldi, P.~Martinon, and J.~Rouot, ``{The Purcell
  Three-link swimmer: some geometric and numerical aspects related to periodic
  optimal controls},'' Oct. 2015, working paper or preprint. [Online].
  Available: \url{https://hal.inria.fr/hal-01143763}
\BIBentrySTDinterwordspacing

\bibitem{Giraldi:2013aa}
L.~Giraldi, P.~Martinon, and M.~Zoppello, ``Controllability and optimal strokes
  for n-link microswimmer,'' in \emph{Decision and Control (CDC), 2013 IEEE
  52nd Annual Conference on}, Dec 2013, pp. 3870--3875.

\bibitem{DeSimone:2012aa}
A.~DeSimone, L.~Heltai, F.~Alouges, and A.~Lefebvre-Lepot, ``Computing optimal
  strokes for low {R}eynolds number swimmers,'' \emph{Natural Locomotion in
  Fluids and on Surfaces}, pp. 177--184, 2012.

\bibitem{Bryson:1979}
\BIBentryALTinterwordspacing
A.~E. Bryson and Y.-C. Ho, \emph{Applied Optimal Control, Optimization,
  Estimation, and Control}.\hskip 1em plus 0.5em minus 0.4em\relax New
  York-London-Sydney-Toronto: John Wiley \& Sons, 1975. [Online]. Available:
  \url{http://dx.doi.org/10.1002/zamm.19790590826}
\BIBentrySTDinterwordspacing

\bibitem{Ostrowski:2000vl}
J.~P. Ostrowski, J.~P. Desai, and V.~Kumar, ``Optimal gait selection for
  nonholonomic locomotion systems,'' \emph{International Journal of Robotics
  Research}, vol.~19, no.~3, pp. 225--237, 2000.

\bibitem{Becker:2003}
L.~Becker, S.~A. Koehler, and H.~A. Stone, ``On self-propulsion of
  micro-machines at low {R}eynolds number: {P}urcell's three-link swimmer,''
  \emph{Journal of Fluid Mechanics}, vol. 490, pp. 15--35, 2003.

\bibitem{Magnus:1954vl}
W.~Magnus, ``On the exponential solution of differential equations for a linear
  operator,'' \emph{Communications on Pure and Applied Mathematics}, vol. VII,
  pp. 649--673, 1954.

\bibitem{Flanders:1973aa}
H.~Flanders, ``Differentiation under the integral sign,'' \emph{The American
  Mathematical Monthly}, vol.~80, no.~6, pp. 615--627, 6 1973.

\bibitem{Kobayashi:1963}
S.~Kobayashi and K.~Nomizu, \emph{Foundations of Differential Geometry}.\hskip
  1em plus 0.5em minus 0.4em\relax Wiley Interscience, 1963, vol.~1.

\end{thebibliography}

\end{document}